\documentclass[12pt]{article}
\usepackage{amsfonts}
\usepackage{amssymb}
\newtheorem{theorem}{Theorem}[section]
\newtheorem{lemma}[theorem]{Lemma}

\newtheorem{corollary}[theorem]{Corollary}

\def\t{\mbox{tr}\,}   
\def\zC{\mathbb C}    
\def\zR{\mathbb R}    
\def\qed{$\hfill\square$}

\title{Young's Inequality in Semifinite von Neumann Algebras}

\author{Douglas R.~Farenick and S.~Mahmoud Manjegani\\
{ } \\
{\em
Department of Mathematics and Statistics,
University of Regina}\\
{\em Regina, Saskatchewan S4S 0A2, Canada} }

\begin{document}
\maketitle


\begin{abstract} This paper formulates Young-type inequalities
for singular values (or $s$-numbers) and traces
in the context
of von Neumann algebras. In particular, it shown that
if $\t(\cdot)$ is a faithful
semifinite normal trace on a semifinite von Neumann algebra $M$ and if
$p$ and $q$ are positive real numbers for which $p^{-1}+q^{-1}=1$,
then, for all positive operators $a,b\in M$,
$\t(|ab|)\le p^{-1}\t(a^p)+ q^{-1}\t(b^q)$,
with
equality holding (in the cases where
$p^{-1}\t(a^p)+ q^{-1}\t(b^q)<\infty$)
if and only if $b^q=a^p$.
\end{abstract}
\maketitle

\section{Introduction}\label{S:introduction}

Young's inequality (see, for example,
page 17 of \cite{hardy--littlewood--polya}) asserts
that if $p$ and $q$ are positive real numbers for which $p^{-1}+q^{-1}=1$,
then
$|\lambda\mu|\le p^{-1}|\lambda|^p+q^{-1}|\mu|^q$, for all complex numbers
$\lambda$ and $\mu$, and equality holds if and only if
$|\mu|^q=|\lambda|^p$.
Several
generalisations of Young's inequality whereby $\lambda$
and $\mu$ are replaced by Hilbert space
operators---or by singular values, norms, or traces
of operators---are known \cite{ando},
\cite{bhatia--parthasarathy}, \cite{erlijman--farenick--zeng}, \cite{kosaki}.
The present paper adds to these results by
formulating new Young-type inequalities in the context
of von Neumann algebras.

In what follows, $N$ shall denote an
arbitrary von Neumann algebra of operators acting on a complex
infinite-dimensional Hilbert space $H$. By $M$ we denote a semifinite
von Neumann algebra and $\t(\cdot)$ is assumed to be a faithful semifinite
normal trace on $M$.
The cone of positive operators in $N$ and the projection lattice
in $N$ are denoted by $N^+$ and $\mathcal P(N)$ respectively.
For any $z\in N$, $|z|$ denotes
$(z^*z)^{1/2}$, the unique positive
square root of $z^*z$.
The notation $e\sim f$, for $e,f\in\mathcal P(N)$, shall indicate
that $e$ and $f$ are Murray--von Neumann equivalent, which is to say that
$e=v^*v$
and $f=vv^*$ for some $v\in N$. If $x\in N$, then $R[x]$ denotes
the range projection for $x$ (that is,
the projection
in $N$ whose range is the closure of the range of $x$.)
The spectral resolution of the
identity of $a\in N^+$ is denoted by $p^a$, yielding
\[a\;=\;\int_{0}^{\infty}s\,dp^a(s)\;.\]

Murray and von Neumann, in their 1936 paper \cite{murray--von neumann},
introduced generalised
singular values, or $s$-numbers,
for operators in semifinite von Neumann algebras. Years later, in the
1980s,
interest in majorisation and operator
inequalities---some of which is chronicled in
\cite{bhatia} and \cite{zhan}---led to careful studies
of singular values by Fack \cite{fack},
Fack and Kosaki \cite{fack--kosaki}, Hiai and
Nakamura \cite{hiai--nakamura},
and others, in the operator algebra context.

The singular values $\mu_z(t)$ of $z\in M$ are
defined for each $t\in\zR_0^+$, where
$\zR_0^+$ is the set of nonnegative
real numbers, by the equation
\[
\mu_z(t)\;=\;\mbox{inf}\,\left\{\|ze\|\,:\,e\in\mathcal P(M),\,\t(1-e)\le t
\right\}\;.
\]
The trace of $|x|$, for $x\in M$, is recovered from the singular values
of $x$ by way of the equation
\[
\t(|x|)\;=\; \int_0^{\infty}\mu_x(t)\,dt\,.
\]
Thus, inequalities between each of the singular values of a pair of
positive operators necessarily imply an inequality between the traces of these
operators. In consequence, one aims to formulate operator inequalities at the
level of singular values, if possible.

The singular
value inequality to be
established in the present paper is
inequality (\ref{E:main inequality 1}) below.
Inequality (\ref{E:main inequality 1})
was first established by Ando \cite{ando}
for finite-dimensional $M$, whereas the cases of equality,
again for $M$ of finite dimension, were analysed
by Hirzallah and Kittaneh \cite{hirzallah--kittaneh}.

\begin{theorem}\label{main result in singular values}
Assume that $M$ is a semifinite von Neumann algebra and that $\t(\cdot)$
is a faithful semifinite normal trace on $M$. Let $p,q\in\zR^+$ satisfy
$p^{-1}+q^{-1}=1$. Then, for every $x,y\in M$,
\begin{equation}\label{E:main inequality 1}
\mu_{|xy^*|}(t)\;\le\;\mu_{p^{-1}|x|^p+q^{-1}|y|^q}(t)\;,
\;\mbox{ for all }\;t\in\zR_0^+\;.
\end{equation}
If $\t(1)<\infty$, then equality holds in the Young
inequality $($\ref{E:main inequality 1}$)$,
for some $x,y\in M$, if
and only if $|y|^q=|x|^p$.
\end{theorem}

Young inequalities in traces can also be formulated, leading to
the following result
for C$^*$-algebras.

\begin{theorem}\label{main result for traces}
Assume that $A$ is a unital C$^*$-algebra and that $\tau$
is a faithful tracial state on $A$. Let $p,q\in\zR^+$ satisfy
$p^{-1}+q^{-1}=1$. Then, for every $a,b\in A^+$,
\begin{equation}\label{E:main inequality 2}
\tau(|ab|)\;\le\;p^{-1}\tau(a^p)\,+\,q^{-1}\tau(b^q)\;.
\end{equation}
Equality holds in the Young
inequality $($\ref{E:main inequality 2}$)$,
for some $a,b\in A^+$, if
and only if $b^q=a^p$.
\end{theorem}

The proofs
of
Theorems \ref{main result in singular values}
and \ref{main result for traces} make extensive
use of various properties of singular values;
these propeties are review below.
Further details can be found in \cite{fack}.

For each $z\in M$, the function $\mu_z:\zR_0^+\rightarrow\zR_0^+$
is nonincreasing and continuous on the right.
Moreover,
$\mu_z=\mu_{z^*}=\mu_{|z|}$ and, consequently, for any $x,y\in M$,
\begin{equation}\label{E:symmetry identity}
\mu_{|xy^*|}\;=\;\mu_{|yx^*|}\;.
\end{equation}
In addition, for every
$w_1,w_2\in M$ and $t\in\zR_0^+$,
\begin{equation}\label{E:contraction inequality}
\mu_{w_1zw_2}(t)\le\|w_1\|\,\|w_2\|\,\mu_z(t)\,.
\end{equation}
The dependence of $\mu_z$ on $z$ is as follows:
if $z_1,z_2\in M$, then, for all $t\in\zR_0^+$,
\begin{equation}\label{E:approximation}
\left|\, \mu_{z_1}(t)\,-\,\mu_{z_2}(t)\,\right|\;\le\;\|z_1\,-\,z_2\|\;.
\end{equation}

The singular values of positive
operators are especially well behaved.
If $h\in M^+$, then
\begin{equation}\label{E:spectral scale 1}
\mu_h(t)\;=\;\min\,\left\{s\in\zR_0^+\,:
\,\t\left(p^h(s,\infty)\right)\le t\right\}\,.
\end{equation}
Alternatively, one can employ a variational principle
to evaluate $\mu_h(t)$:
\[
\mu_h(t)=
\mbox{inf}\left\{\mbox{sup}\{\langle h\xi,\xi\rangle:
\xi\in\mbox{ran}\,e,\,\|\xi\|=1\}:
e\in\mathcal P(M),\,\t(1-e)\le t\right\}\,.
\]
Finally, we shall also require the following
continuous functional calculus (Proposition 1.6 of
\cite{fack}).
If $\psi:\zR_0^+\rightarrow\zR_0^+$ is an increasing continuous function
such that $\psi(0)=0$, then
\begin{equation}\label{E:functional calculus}
\mu_{\psi(h)}(t)\;=\;\psi\left(\,\mu_h(t)\,\right)\,,\,\mbox{ for all }\,t\in\zR_0^+\,.
\end{equation}

Throughout, the notation $a\le b$, for hermitian operators
$a,b\in N$, refers to the L\"owner partial order, namely
$a\le b$ if and only if
$\langle a\xi,\xi\rangle\le\langle b\xi,\xi\rangle$ for all $\xi\in H$.


\section{Inequalities}
\label{S:inequalties}

Ando proved in \cite{ando} that, for
operators on finite-dimensional Hilbert spaces,
Young's inequality holds at
the level of singular values.
This Young-type inequality for singular values
was later extended to compact operators in
\cite{erlijman--farenick--zeng}. Theorem \ref{singular values}
below represents the most general
form of Ando's original result.
However, the proof of Theorem \ref{singular values} and
other results herein
rest upon a core result concerning a compressed
form of Young's inequality. This ``compression lemma" was also first
established in finite dimensions by Ando \cite{ando}, but it holds
in arbitrary von Neumann algebras as well \cite{erlijman--farenick--zeng}.

\begin{lemma}[Compression Lemma] \label{compression lemma}
Assume that $p\in(1,2]$, $q=(1-p^{-1})^{-1}$, $a,b\in N^+$, and
$b$ is invertible.
If
$f_s=R[b^{-1}p^{|ab|}(\,(s,\infty)\,)]$,
for
$s\in\zR_0^+$,
then
\[
sf_s\;\le\;f_s\left(p^{-1}a^p\,+\,q^{-1}b^q\right)f_s\quad\mbox{ and }\quad
f_s\sim p^{|ab|}(\,(s,\infty)\,)\;.
\]
\end{lemma}

\noindent{\em Proof.} Except for the claim that
$f_s\sim p^{|ab|}(\,(s,\infty)\,)$, the rest of the
lemma is precisely
Proposition 2.3 of \cite{erlijman--farenick--zeng}.
To prove that $f_s\sim p^{|ab|}(\,(s,\infty)\,)$, it is enough to
prove the following general proposition:
if $e,f\in\mathcal P(N)$, $b\in N^+$ is invertible, and
$f=R[b^{-1}e]$, then $e\sim f$.
To this end, if $b^{-1}e=v|b^{-1}e|$ is the polar decomposition
of $b^{-1}$, then
$v^*v=R[eb^{-1}]=e$ and
$vv^*=R[b^{-1}e]=f$, whence $e\sim f$.
\qed\medskip

\begin{lemma}\label{spectral scale of a projection} If $f\in\mathcal P(M)$,
then $\mu_f(t)=1$ for all $t<\t(f)$
and $\mu_f(t)=0$ for all $t\ge\t(f)$.
\end{lemma}

\noindent{\em Proof.}
If $t\ge\t(f)$, then set $e=1-f$ to obtain $\t(1-e)\le t$
and
$0=\|f(1-f)\|=\|fe\|\ge\mu_f(t)\ge0$,
which shows that $\mu_f(t)=0$.

Assume now that $t<\t(f)$.
If $e\in\mathcal P(M)$ satisfies
$\t(1-e)\le t$, then $\t(1-e)<\t(f)$.
By Kaplansky's Lemma,
\[
e-(e\wedge f)\;\sim\;(e\vee f)-f\;\le\;1-f\,,
\]
and so if it were true that $e\wedge f=0$, then we
would have $\t(e)\le\t(1-f)$, or equivalently
$\t(1-e)\ge\t(f)$, in contradiction to
$\t(1-e)<\t(f)$. Thus, it must be that
$e\wedge f\not=0$, and so $\|fe\|=1$.
This proves that $\|fe\|=1$ for all $e\in\mathcal P(M)$
that satisfy $\t(1-e)\le t$;
hence,
$\mu_f(t)=1$.
\qed\medskip

The general form of Ando's theorem \cite{ando} can now be established.

\begin{theorem}[Young's Inequality in Singular Values]\label{singular values} If
$p$ and $q$ are positive real numbers for which $p^{-1}+q^{-1}=1$,
and if $x,y\in M$ and $t\in\zR_0^+$, then
\begin{equation}\label{E:young's inequality in singular values}
\mu_{|xy^*|}(t)\;\le\; \mu_{p^{-1}|x|^p\,+\,q^{-1}|y|^q}(t)\;.
\end{equation}
\end{theorem}

\noindent{\em Proof.}
We begin by
showing that the proof can be reduced to the case of
positive operators.

If $y=w|y|$ is the polar decomposition of $y$ in $M$,
then the proof of Proposition 4.1 in \cite{erlijman--farenick--zeng}
demonstrates that
\begin{equation}\label{E:polar decomposition identity}
|xy^*|=w
\left|\phantom{\int}\!\!\!\!\!
|x||y| \phantom{\int}\!\!\!\!\! \right|
w^*\;.
\end{equation}
By $\|w\|\le1$ and
property (\ref{E:contraction inequality}) of $\mu_z$, we have that
\begin{equation}\label{E:pass to positives}
 \mu_{|xy^*|}(t)\le \mu_{\left|\phantom{\int}\!\!\!\!\!
\,|x||y| \,\phantom{\int}\!\!\!\!\! \right|}(t)\;.
\end{equation}
Thus,
in setting $a=|x|$ and $b=|y|$, it is sufficient to prove that
\begin{equation}\label{E:reduction to positives}
\mu_{|ab|}(t)\;\le\; \mu_{p^{-1}a^p\,+\,q^{-1}b^q}(t)\;.
\end{equation}
Therefore, we shall prove that inequality (\ref{E:reduction to positives})
holds for all $a,b\in M^+$.

So, we assume henceforth that $a,b\in M^+$. We assume, further, that
$p\in(1,2]$ and that $b\in M^+$ is invertible. The assumption on $p$
entails no loss of generality
because if inequality
(\ref{E:reduction to positives}) holds for $1<p\le 2$, then in
cases where $p>2$ the conjugate $q$ satisfies $q<2$, and so
\[ \mu_{|ab|}(t)\;=\;\mu_{|ba|}(t)\;\le\;\mu_{q^{-1}b^q+p^{-1}a^p}(t)\;,\]
where the equality
$\mu_{|ab|}=\mu_{|ba|}$ is obtained from property (\ref{E:symmetry identity})
of the function $\mu_z$.
The assumption that $b\in M^+$ be invertible is also no loss in
generality, for if $b$ is not invertible, then
consider $b_{\varepsilon}=b+\varepsilon1$, an invertible element
for which $\|b_{\varepsilon}-b\|\rightarrow 0$ as $\varepsilon\rightarrow 0^+$.
Property (\ref{E:approximation}) of the function $\mu_z$ implies that
\[ \lim_{\varepsilon\rightarrow 0^+}\mu_{|ab_{\varepsilon}|}(t)\;=\;\mu_{|ab|}(t)\;,\]
for every $t\in\zR_0^+$. Thus, if inequality
(\ref{E:reduction to positives}) holds for invertible elements, then for every $t$,
\[\mu_{|ab|}(t)\;=\; \lim_{\varepsilon\rightarrow 0^+}\mu_{|ab_{\varepsilon}|}(t)
\;\le\;\lim_{\varepsilon\rightarrow 0^+}
\mu_{p^{-1}a^p\,+\,q^{-1}b_{\varepsilon}^q}(t)
\;=\; \mu_{p^{-1}a^p\,+\,q^{-1}b^q}(t)\;.\]
Hence, it is enough to already assume that $b\in M^+$ is invertible.

Invoke Lemma \ref{compression lemma} to obtain, for
$s\in\zR_0^+$ and
$f_s=R[b^{-1}p^{|ab|}(\,(s,\infty)\,)]$,
\begin{equation}\label{E: compression inequality}
sf_s\;\le\;f_s\left(p^{-1}a^p\,+\,q^{-1}b^q\right)f_s
\;
\end{equation}
and $f_s\sim p^{|ab|}(\,(s,\infty)\,)$.
Consequently, $\t(f_s)=\t(\,p^{|ab|}(\,(s,\infty)\,)\,)$.

Now
fix
$t\in\zR_0^+$
and let $\zeta=\mu_{|ab|}(t)$.
If $\zeta=0$, then inequality (\ref{E:reduction to positives})
holds trivially. Therefore, assume that $\zeta>0$.
Suppose that $\varepsilon>0$
satisfies $(\zeta-\varepsilon)>0$.
By (\ref{E:spectral scale 1}),
\[ \zeta\;=\;\min\,\left\{s\in\zR\,:
\,\t\left(p^{|ab|}(s,\infty)\right)\le t\right\}\,. \]
With $s=\zeta$, Lemma
\ref{compression lemma}
yields $f_{\zeta}\sim
p^{|ab|}(\,(\zeta,\infty)\,)$. Hence,
\[ \t(f_{\zeta})\;=\;\t\left(p^{|ab|}(\zeta,\infty)\right)
\;\le\;t\,.\]
Thus,
\[ 0\le(\zeta-\varepsilon)<\zeta \quad\Longrightarrow\quad
\t(f_{\zeta-\varepsilon})\,>\,t\;.\]
Now replace in $s$ with $\zeta-\varepsilon$
in inequality (\ref{E: compression inequality}) to
obtain the inequality
\begin{equation}\label{E:epsilon compression inequality 1}
(\zeta-\varepsilon)f_{\zeta-\varepsilon}\;\le\;f_{\zeta-\varepsilon}(
p^{-1}a^p\,+\,q^{-1}b^q)f_{\zeta-\varepsilon}\;.
\end{equation}
By (\ref{E:contraction inequality}),
inequality (\ref{E:epsilon compression inequality 1})
yields
\begin{equation}\label{E: epsilon compression inequality}
(\zeta-\varepsilon)\mu_{f_{\zeta-\varepsilon}}(t)\;\le\;\mu_{p^{-1}a^p+q^{-1}b^q}(t)
\quad\mbox{ for all }t\in\zR_0^+\,.
\end{equation}
By Lemma \ref{spectral scale of a projection}, $\t(f_{\zeta-\varepsilon})>t$
implies that $\mu_{f_{\zeta-\varepsilon}}(t)=1$. Therefore,
inequality (\ref{E: epsilon compression inequality}) can be rewritten as
\[
\zeta\,\le\,\mu_{p^{-1}a^p+q^{-1}b^q}(t)\,+\,\varepsilon\;.
\]
Because $\mu_{|ab|}(t)\,=\,\zeta$ and because
the inequality above is true for every $\varepsilon>0$,
$\mu_{|ab|}(t)\le\mu_{p^{-1}a^p+q^{-1}b^q}(t)$, which completes the
proof.
\qed\medskip

Theorem \ref{singular values}
does not hold, in general, if $\mu_{|xy^*|}(t)$
is replaced by
$\mu_{|xy|}(t)$ on the left hand side. A counterexample
in $2\times 2$ matrices can be found on p.~263 of \cite{bhatia}.

As noted in the introduction, Young's inequality in singular
values automatically leads to a Young-type inequality for
traces.

\begin{corollary}[A Tracial Young Inequality] \label{traces} If
$p$ and $q$ are positive real numbers for which $p^{-1}+q^{-1}=1$,
then, for all $x,y\in M$,
\begin{equation}\label{E:in traces}
\t(|x^{\dagger}y^{\dagger}|)
\;\le\; p^{-1}\t(|x|^p)\,+\,q^{-1}\t(|y|^q)\;,
\end{equation}
where $x^{\dagger}\in\{x,x^*, |x|,|x^*|\}$, $y^{\dagger}\in\{y,y^*,|y|,|y^*|\}$.
\end{corollary}

\noindent{\em Proof.} Theorem \ref{singular values} and
the integral representation of traces
imply that
\[
\begin{array}{rcl}
\t(|xy^*|)\;&=&\;\displaystyle\int_0^{\infty}\mu_{|xy^*|}(t)\,dt \\
\;&\le&\;
\displaystyle\int_0^{\infty}\mu_{p^{-1}|x|^p+q^{-1}|y|^q}(t)\,dt \\
\;&=&\;
\t(p^{-1}|x|^p+q^{-1}|y|^q)\;.
\end{array}
\]
Thus, what remains to be shown is that
the right hand side of the inequality does not
change if $x$ and $y$ are replaced by $x^{\dagger}$ and
$y^{\dagger}$. That is, we need only show that
$\t(|x^{\dagger}|^p)=\t(|x|^p)$ and that $\t(|y^{\dagger}|^q)
=\t(|y|^q)$. It is enough to consider the case of $x$. The
identity $\t(x^*x)=\t(xx^*)$ implies that $\t((x^*x)^k)=\t((xx^*)^k)$
for all positive integers $k$. Thus, by functional calculus,
$\t(\psi(|x|))=\t(\psi(|x^*|))$,
for all continuous functions $\psi:\zR_0^+\rightarrow\zR_0^+$---and
in particular for
$\psi(t)=t^p$.
\qed\medskip

An elegant and
far-reaching theory of majorisation in von Neumann algebras
was developed by Hiai in \cite{hiai}. If $x,y\in M$,
then we write $\mu(x)\prec_w\mu(y)$ to denote that the
singular values of $x$ are weakly majorised by the
singular values of $y$. That is,
\[
\mu(x)\prec_w\mu(y)\,\mbox{ if and only if }\,
\int_0^s\mu_x(t)\,dt\;\le\;\int_0^s\mu_y(t)\,dt\;\mbox{ for all }\,
s\in\zR_0^+\;.
\]
The following result
(Proposition 4.3 of \cite{fack}) concerning weak majorisation is
fundamental: $\mu(xy)\prec_w\mu(x)\mu(y)$. That is,
\begin{equation}\label{E:submajorisation}
\int_0^{s}\mu_{xy}(t)\,dt\;\le\; \int_0^{s}\mu_x(t)\mu_y(t)\,dt\;
\mbox{ for all }\,s\in\zR_0^+\,.
\end{equation}
Of course, Theorem \ref{singular values} implies
that
$\mu(|xy^*|)\prec_w\mu(p^{-1}|x|^p+q^{-1}|y|^q)$,
for all $x,y\in M$.

In addition to weak majorisation, it is useful to consider
the spectral pre-order, which can formulated in arbitrary
von Neumann algebras.
For $a,b\in N^+$, we write $a \prec_{sp} b$
to indicate
that $p^a(s,\infty)$ is Murray-von Neumann equivalent to a subprojection
of $p^b(s,\infty)$, for every $s\in\zR_0^+$.

We now
apply a little of the theory
of majorisation
to the singular value Young inequality to obtain
Theorem \ref{hiai} below, which shows that, in
finite factors,
Young's inequality in singular
values implies a Young inequality in the spectral
pre-order and a Young inequality
in the
L\"owner partial order---after
a correction by a doubly stochastic map.
(A positive linear map $\Phi:M\rightarrow M$ is doubly stochastic if
$\Phi(1)=1$ and $\t\left(\Phi(h)\right)=\t(h)$ for all $h\in M^+$.)

\begin{theorem}\label{hiai} Assume that $p$ and $q$
are positive real numbers for which $p^{-1}+q^{-1}=1$,
and let $x,y\in M$, where $M$ is a (semifinite) factor.
\begin{enumerate}
\item $|xy^*|\prec_{sp}p^{-1}|x|^p+q^{-1}|y|^q$, and
\item if $M$ is finite, then
$|xy^*|\le p^{-1}\Phi(|x|^p)+q^{-1}\Phi(|y|^q)$ for some
doubly stochastic positive linear map $\Phi:M\rightarrow M$.
\end{enumerate}
\end{theorem}

\noindent{\em Proof.} Because $M$ is a factor,
if $e,f\in\mathcal P(M)$, then $e$ is equivalent to a
subprojection of $f$ or vice versa. If we consider this
fact with the spectral projections
\[
p^{|xy^*|}(s,\infty)\quad\mbox{ and }\quad
p^{p^{-1}|x|^p+q^{-1}|y|^q}(s,\infty)\,,
\]
for all $s\in\zR_0^+$, then
the Young inequality
$\mu_{|xy^*|}(t)\le\ \mu_{p^{-1}|x|^p+q^{-1}|y|^q}(t)$,
(Theorem \ref{singular values}) yields
$|xy^*|\prec_{sp}p^{-1}|x|^p+q^{-1}|y|^q$, which proves
the first assertion.

If $M$ is finite, then
Young's inequality $\mu_{|xy^*|}(t)\le\ \mu_{p^{-1}|x|^p+q^{-1}|y|^q}(t)$,
for all $t\in\zR_0^+$, and Theorem 4.7 of \cite{hiai}
imply that $|xy^*|=\Psi(p^{-1}|x|^p+q^{-1}|y|^q)$,
for some positive linear map $\Psi:M\rightarrow M$ for
which $\t(\Psi(a))\le\t(a)$, for all $a\in M^+$.
Proposition 4.3 of \cite{hiai}
completes the argument:
because $M$ is finite, there is a
doubly stochastic positive linear map $\Phi:M\rightarrow M$
such that $\Psi(a)\le\Phi(a)$, for all $a\in M^+$.
\qed\medskip

It would be
interesting to know whether the doubly stochastic positive linear
map $\Phi:M\rightarrow M$
could in fact be chosen to be an
automorphism. (This is the case if $M$ is a factor of type
I${}_n$, as the spectral theorem and Theorem \ref{singular values}
demonstrate.)

In the case $p=q=2$, Young's inequality
is the arithmetic--geometric
mean inequality, which, if $M$ were taken to be $\zC$, would be in the formulation
$\alpha\beta\le\frac{1}{2}(\alpha^2+\beta^2)$, for positive real numbers
$\alpha,\beta$. Sometimes, however, one wants the arithmetic--geometric
mean inequality in its more traditional (equivalent) form: $\sqrt{\alpha\beta}\le
\frac{1}{2}(\alpha+\beta)$. Such a formulation extends to noncommutative $M$
in Theorem \ref{agm} below as weak majorisation and as
a tracial inequality .
(See, also, \cite{bhatia-jot} for a formulation in unitarily-invariant
norms.)

\begin{theorem}[Arithmetic--Geometric Mean Inequality]\label{agm} If
$a,b\in M^+$, then
\[
\mu(|ab|^{1/2}) \;\prec_w\;
\frac{1}{2}\left(\,\mu(a)+\mu(b)\,\right)
\]
and
\[
\t(|ab|^{1/2})\;\le\;\left(\t(a)\t(b)\right)^{1/2}\;\le\;
\frac{1}{2}(\,\t(a)+\t(b)\,)\;.
\]
\end{theorem}

\noindent{\em Proof.} For any $h\in M^+$, let $\Lambda_h:(0,\t(1))\rightarrow\zR^+$
denote the function
\[
\Lambda_h(s)\;=\;\mbox{exp}\left(\int_0^s\mbox{log}\,\mu_h(t)\,dt\right)\,.
\]
Because $\mu_{h^{1/2}}(t)=\sqrt{\mu_h(t)}$ for all $t\ge0$,
by the functional calculus
(\ref{E:functional calculus}), it follows that
$\Lambda_{h^{1/2}}(s)=\sqrt{\Lambda_h(s)}$ for all $s\in(0,\t(1))$.
Furthermore, Theorem 2.3 of \cite{fack} indicates that
$\Lambda_{|ab|}(s)\le\Lambda_a(s)\Lambda_b(s)$ and so
$\Lambda_{|ab|^{1/2}}(s)\le \Lambda_{a^{1/2}}(s)\Lambda_{b^{1/2}}(s)$.
Using this last inequality and the equations
\[
\begin{array}{rcl}
\Lambda_{a^{1/2}}(s)\Lambda_{b^{1/2}}(s)\;&=&\;
\mbox{exp}\left(
\displaystyle\int_0^s(\mbox{log}\,\mu_{a^{1/2}}(t)\;+\;\mbox{log}\,\mu_{b^{1/2}}(t))dt
\right)\\
\;&=&\;
\mbox{exp}\left(
\displaystyle\int_0^s\mbox{log}\,\sqrt{\mu_a(t)\mu_b(t)}\,dt
\right)\,,
\end{array}
\]
we see that, for all $s\in(0,\t(1))$,
\begin{equation}\label{E:integral inequality 1}
\int_0^s\mbox{log}\,\mu_{|ab|^{1/2}}(t)\,dt\;\le\;
\int_0^s\mbox{log}\,\sqrt{\mu_a(t)\mu_b(t)}\,dt\;.
\end{equation}
Each of the integrands in (\ref{E:integral inequality 1})
above is nonincreasing, and so if $f$ is any
increasing convex function, then
for all $s\in(0,\t(1))$ (by majorisation theory
\cite{hardy--littlewood--polya}),
\begin{equation}\label{E:integral inequality 2}
\int_0^sf\left(\mbox{log}\,\mu_{|ab|^{1/2}}(t)\right)\,dt\;\le\;
\int_0^sf\left(\mbox{log}\,\sqrt{\mu_a(t)\mu_b(t)}\right)\,dt\;.
\end{equation}
In particular, inequality (\ref{E:integral inequality 2}) holds
for the function $f(r)=e^r$. Hence,
\[
\begin{array}{rcl}
\displaystyle\int_0^{s}\mu_{|ab|^{1/2}}(t)\,dt\;&\le&\;
\displaystyle\int_0^{s}\sqrt{\mu_a(t)\mu_b(t)}\,dt\\
\;&\le&\;
\left(\displaystyle\int_0^{s}\mu_a(t)\,dt\right)^{\frac{1}{2}}
\left(\displaystyle\int_0^{s}\mu_b(t)\,dt\right)^{\frac{1}{2}}\\
\;&\le&\;
\frac{1}{2}
\displaystyle\int_0^{s}\mu_a(t)\,dt\;+\;
\frac{1}{2}\displaystyle\int_0^{s}\mu_b(t)\,dt\,.
\end{array}
\]
(The second inequality above
is the Cauchy--Schwarz inequality.)
\qed\medskip

Theorem \ref{agm} above is proved with considerable less effort
than what
Theorem \ref{singular values} requires. Some other related
tracial inequalities are just as readily established; for example,
the Fenchel--Young inequality.
To describe what is involved,
assume that $F:\zR_0^+\rightarrow\zR_0^+$ is a convex function and that
$F^*:\Gamma_F\rightarrow\zR_0^+$ is the Fenchel conjugate
\cite{fenchel} of $F$. That is, $\Gamma_F$ is a convex subset of $\zR$
and
\[ F^*(r)\;=\;\sup_{t\in\zR_0^+}\left(rt\,-\,F(t)\right)\;.\]
The Fenchel--Young inequality is $\alpha\beta\le F(\alpha)+F^*(\beta)$,
for all $\alpha\in\zR_0^+$,
$\beta\in\Gamma_F$. (If $p>1$ and $F(t)=p^{-1}t^p$, then the
Fenchel conjugate is $F^*(s)=q^{-1}s^q$, where $p^{-1}+q^{-1}=1$,
and $\Gamma_F=\zR_0^+$. Thus,
the Fenchel--Young inequality implies the Young inequality
under study here.)
If $a,b\in M^+$, then the weak majorisation relation
$\mu(|ab|)\prec_w\mu(a)\mu(b)$
suggests that one can
apply the Fenchel--Young inequality pointwise to
the products $\mu_a(t)\mu_b(t)$
to obtain a Fenchel--Young inequality in traces:
\[
\t(|ab|)\;\le\;\t(F(a))+\t(F^*(b))\,,
\]
for all $a,b\in M^+$ for which $\zR_0^+\cap\Gamma_F$
contains the spectrum of $b$.

\section{Cases of Equality}
\label{S:cases of equality}

There are few results about cases of equality in
operator inequalities; however, Young's inequality is
somewhat of an exception.
A characterisation of equality in the singular value
Young inequality (Theorem \ref{singular values}) was given
for
finite-dimensional $M$ in
\cite{hirzallah--kittaneh},
whereas a
characterisation of equality (for elements of finite
trace) in the tracial
Young inequality (Corollary \ref{traces})
was given in \cite{argerami--farenick} for the
I${}_{\infty}$-factor $B(H)$. Theorems \ref{equality in traces}
and
\ref{equality in singular values}
below add to these results by characterising the cases of
equality in the Young inequalities when $M$ is
an arbitrary semifinite von Neumann algebra.

\begin{theorem}[Equality in Traces]
\label{equality in traces}
Let
$p$ and $q$ be positive real numbers such that
$p^{-1}+q^{-1}=1$, and assume that $a,b\in M^+$
satisfy $\t(a)<\infty$ and $\t(b)<\infty$.
Then
\begin{equation}\label{E:equality in traces}
\t(|ab|)\;=\;p^{-1}\t(a^p)\,+\,q^{-1}\t(b^q)
\end{equation}
if and only if $b^q=a^p$.
\end{theorem}

\noindent{\em Proof.}
It is clear that equation (\ref{E:equality in traces}) holds if $b^q=a^p$, and
so we focus on the converse.

First, observe that if $r\in\zR_0^+$ and $r>1$, and if
$h\in M^+$ has finite trace, then
so does $h^r$.
Indeed,
\[
\t(h)\;=\;\int_0^{\infty}\mu_h(t)\,dt\;<\;\infty\,
\]
implies that $\mu_h(t)\rightarrow 0^+$ as $t\rightarrow\infty$.
Thus, for sufficiently large $t$,
$[\mu_h(t)]^r\le\mu_h(t)$. Hence, using
$\mu_{h^r}(t)=[\mu_h(t)]^r$ by (\ref{E:functional calculus}),
$\t(h^r)<\infty$.

Assume that  equation (\ref{E:equality in traces}) holds for some
$a,b\in M^+$.
Inequality (\ref{E:submajorisation}) states
that
\[
\t(|ab|)\;=\;
\int_0^{\infty}\mu_{ab}(t)\,dt\;\le\; \int_0^{\infty}\mu_a(t)\mu_b(t)\,dt\;.
\]
Thus, using that equality (\ref{E:equality in traces}) holds,
\[
\begin{array}{rcl}
p^{-1}\t(a^p)\,+\,q^{-1}\t(b^q)
\;&=&\;
\displaystyle\int_0^{\infty} \left(p^{-1}\mu_{a^p}(t)\,+
\,q^{-1}\mu_{b^q}(t)\right)\,dt \\
\;&\le&\; \displaystyle\int_0^{\infty} \mu_a(t)\mu_b(t)\,dt \;.
\end{array}
\]
By Young's inequality,
\begin{equation}\label{E:Young in mu}
\mu_a(t)\mu_b(t)\;\le\;p^{-1}[\mu_a(t)]^p\,+\,q^{-1}[\mu_b(t)]^q
\;=\;p^{-1}\mu_{a^p}(t)\,+\,q^{-1}\mu_{b^q}(t)
\end{equation}
for every $t\in\zR_0^+$. Therefore,
\[
\int_0^{\infty} \mu_a(t)\mu_b(t)\,dt
\;=\;
\int_0^{\infty}
\left(p^{-1}\mu_{a^p}(t)+q^{-1}\mu_{b^q}(t)\right)\,dt\;.
\]
This shows, when coupled with (\ref{E:Young in mu}), that
\[ \mu_a(t)\mu_b(t)\;=\; p^{-1}\mu_{a^p}(t)\,+\,q^{-1}\mu_{b^q}(t)\,
\]
for almost all $t\in\zR_0^+$. However, as the nonincreasing
functions $\mu_z$, for $z\in M$, are right
continuous,
$\mu_a(t)\mu_b(t)=p^{-1}\mu_{a^p}(t)+q^{-1}\mu_{b^q}(t)$ for all
$t\in\zR_0^+$.
But these are cases of equality in Young's inequality,
and so
$[\mu_b(t)]^q=[\mu_a(t)]^p$
for all $t\in\zR_0^+$. That is,
by again using (\ref{E:functional calculus}),
\[
\mu_{b^q}(t)\;=\;\mu_{a^p}(t)\;\mbox{ for all }\;t\in\zR_0^+\,.
\]
With $t=0$, this equation
implies that
\begin{equation}\label{E:equality in norm}
\|b\|\;=\;\|a\|^{\frac{p}{q}}\;.
\end{equation}
Integration over all $t$ yields
\begin{equation}\label{E:traces are equal}
\t(a^p)\;=\;\int_0^{\infty}\mu_{a^p}(t)\,dt\;=\;
\int_0^{\infty}\mu_{b^q}(t)\,dt\;=\;\t(b^q)\;.
\end{equation}

A similar argument to the one above shows that
$\mu_{|ab|}(t)=\mu_{p^{-1}a^p+q^{-1}b^q}(t)$ for all
$t\in\zR_0^+$.
The reasons for this are:
$\mu_{|ab|}(t)\le\mu_{p^{-1}a^p+q^{-1}b^q}(t)$ for all $t\in\zR_0^+$
and
\[
\begin{array}{rcl}
\displaystyle\int_0^{\infty}\mu_{|ab|}(t)\,dt\;&=&\;\t(|ab|)\;=\;
p^{-1}\t(a^p)\,+\,q^{-1}\t(b^q) \\
\;&=&\;\displaystyle\int_0^{\infty}
\mu_{p^{-1}a^p+q^{-1}b^q}(t)\,dt\;,
\end{array}
\]
whence
$\mu_{|ab|}(t)=\mu_{p^{-1}a^p+q^{-1}b^q}(t)$ for almost all $t\in\zR_0^+$.

Hence, we have thus far proved that, for every $t\in\zR_0^+$,
\begin{equation}\label{E:equality in s values}
\mu_{|ab|}(t)\;=\;\mu_{p^{-1}a^p+q^{-1}b^q}(t)\;=\;
p^{-1}\mu_{a^p}(t)+q^{-1}\mu_{b^q}(t)\,.
\end{equation}

The remainder of the proof carried out
in cases. The method of the first case, in
particular, is inspired by an approach of Hirzallah and Kittaneh
\cite{hirzallah--kittaneh}.

\vskip 4 pt \noindent{\em Case 1}: ($p=q=2$)
Assume that $p=q=2$. Equation (\ref{E:equality in s values}) becomes
$\mu_{|ab|}(t)=\mu_{\frac{1}{2}(a^2+b^2)}(t)$, for all $t\in\zR_0^+$,
implying that
$\t(\frac{1}{4}(a^2+b^2)^2)=\t(|ab|^2)<\infty$ by the functional
calculus (\ref{E:functional calculus}).
The equation
\[
\t\left(\frac{1}{4}(a^2+b^2)^2-\frac{1}{4}(a^2-b^2)^2-a^2b^2\right)\,=\,0
\]
is readily verified by expanding the left hand
side and it
yields
\[
\frac{1}{4}\t\left((a^2-b^2)^2\right)+\t(a^2b^2)
\,=\,
\t\left(\frac{1}{4}(a^2+b^2)^2\right)
\,=\,\t(|ab|^2)\,.
\]
As $\t(a^2b^2)=\t(|ab|^2)$, the equation above
holds only if $\t((a^2-b^2)^2)=0$.
By the faithfulness of the trace and the uniqueness of positive square
roots, it follows that $b=a$.

\vskip 4 pt\noindent{\em Case 2}: ($p<q$)
Assume that $p<q$; then necessarily $q>2>p$.
We first aim to show that
\begin{equation}\label{E:inequality 1}
\mu_{|ab|}(t)\;\le\;  \mu_{a^{p/2}b^{q/2}}(t)
\|a\|^{1-\frac{p}{2}}  \|b\|^{1-\frac{q}{2}}\;,
\;\mbox{ for all }\;t\in\zR_0^+\;.
\end{equation}
As $1-\frac{q}{2}<0$, the operator $b^{1-\frac{q}{2}}$ exists only
if $b$ is invertible.
Therefore, to prove (\ref{E:inequality 1}) we shall assume that $b$ is invertible.
This assumption entails no loss in generality, for
if $b$ were not invertible, then we could replace $b$ by $b_{\varepsilon}=b+\varepsilon 1$,
which is invertible and which satisfies
$\mu_{zb_{\varepsilon}}(t)\rightarrow
\mu_{zb}(t)$, for every $z\in M$ and $t\in\zR_0^+$, as $\varepsilon\rightarrow 0^+$.
Thus, (\ref{E:inequality 1}) is achieved for noninvertible $b$ as a limiting
case of (\ref{E:inequality 1}) using invertible $b_{\varepsilon}$.
Factor $ab$ as $ab=a^{1-\frac{p}{2}}a^{\frac{p}{2}}b^{\frac{q}{2}}b^{1-\frac{q}{2}}$.
Inequality (\ref{E:contraction inequality}) shows, therefore, that
\[
\mu_{|ab|}(t)\;=\;\mu_{ab}(t)
\;=\;\mu_{ a^{1-\frac{p}{2}}a^{\frac{p}{2}}b^{\frac{q}{2}}b^{1-\frac{q}{2}} }(t)\;\le\;
\|a\|^{1-\frac{p}{2}}\mu_{a^{p/2}b^{q/2}}(t)\|b\|^{1-\frac{q}{2}}\;,
\]
thereby proving (\ref{E:inequality 1}).

Now let $a_1=a^{p/2}$ and $b_1=b^{q/2}$. Because
$\t(a^p)=\t(b^q)$ (equation (\ref{E:traces are equal})), we have that
\begin{equation}\label{E:inequality 2}
\begin{array}{rcl}
\t(|a_1b_1|)\;&\le&\;\displaystyle\frac{1}{2}\t(a_1^2)\,+\,
\displaystyle\frac{1}{2}\t(b_1^2)  \\
& & \\
\;&=&\;
\displaystyle\frac{1}{2}\t(a^p)\,+\,\displaystyle\frac{1}{2}\t(b^q) \\
& & \\
\;&=&\;p^{-1}\t(a^p)\,+\,q^{-1}\t(b^q) \\
& & \\
\;&=&\;
\t(|ab|)\;.
\end{array}
\end{equation}
On the other hand, inequality (\ref{E:inequality 1}) and equation
(\ref{E:equality in norm}) yield
\[
\mu_{|ab|}(t)\;\le\;  \mu_{a^{p/2}b^{q/2}}(t)
\|a\|^{1-\frac{p}{2}}  \|b\|^{1-\frac{q}{2}}\;=\;
\mu_{|a_1b_1|}(t)\|a\|^{1+\frac{p}{q}-p}\;=\;
\mu_{|a_1b_1|}(t)\;.
\]
Thus, upon integrating, we have $\t(|ab|)\le\t(|a_1b_1|)$,
making
inequality (\ref{E:inequality 2}) an equality:
\[
\t(|a_1b_1|)=\frac{1}{2}(\t(a_1^2)+\t(b_1^2))\,.
\]
Hence, from
what was proved in Case 1, we conclude
that $b_1=a_1$ and, as a result, that $b^q=a^p$, as desired.

\vskip 4 pt\noindent{\em Case 3}: ($p>q$)
This case is handled in analogous way to the case $p<q$.
\qed\medskip

The characterisation of equality in the tracial Young
inequality leads to a similar result for singular
values.

\begin{theorem}[Equality in Singular Values]
\label{equality in singular values}
Let $p$ and $q$ be positive real numbers for
which $p^{-1}+q^{-1}=1$. Assume
that $x,y\in M$ satisfy $\t(|x|)<\infty$
and $\t(|y|)<\infty$.
Then
\[ \mu_{|xy^*|}(t)\;=\;\mu_{p^{-1}|x|^p+q^{-1}|y|^q}(t)\;,\;
\mbox{ for all }\;t\in\zR_0^+\,,
\]
if and only if $|y|^q=|x|^p$.
\end{theorem}

\noindent{\em Proof.} Integration leads to
\begin{equation}\label{E:integration equation}
\begin{array}{rcl}
\t(|xy^*|)\;&=&\;\displaystyle\int_0^{\infty}
\mu_{|xy^*|}(t)\,dt \\
\;&=&\;\displaystyle\int_0^{\infty}
\mu_{p^{-1}|x|^p+q^{-1}|y|^q}(t)\,dt  \\
\;&=&\;
p^{-1}\t(|x|^p)\,+\,q^{-1}\t(|y|^q)\;.
\end{array}
\end{equation}
Let $y=w|y|$ be the polar decomposition of $y$.
Then
\[
\begin{array}{rcl}
\mbox{(i)}\;&\;
   w
\left|\phantom{\int}\!\!\!\!\!
|x||y| \phantom{\int}\!\!\!\!\! \right|
w^*  \;&=\; |xy^*| \\
\mbox{(ii)}\;&\;
w^*w\left|\phantom{\int}\!\!\!\!\!
|x||y| \phantom{\int}\!\!\!\!\! \right|
\;&=\;\left|\phantom{\int}\!\!\!\!\!
|x||y| \phantom{\int}\!\!\!\!\! \right| \\
\mbox{(iii)}\;&\;
\t\left( \left|\phantom{\int}\!\!\!\!\!
|x||y| \phantom{\int}\!\!\!\!\! \right|
 \right)
\;&=\;p^{-1}\t(|x|^p)\,+\,q^{-1}\t(|y|^q)\;.
\end{array}
\]
Equation (i) above
is proved in Proposition 4.1 in \cite{erlijman--farenick--zeng}; equation (ii)
is from the polar decomposition;
and
equation (iii) follows from (i), (ii), and (\ref{E:integration equation}).
Hence, by Theorem \ref{equality in traces},
$|y|^q=|x|^p$.
\qed\medskip

\section{Tracial Inequalities in C$^*$-algebras}

By way of
standard arguments with
the Gelfand--Naimark--Segal representation of states,
the tracial Young inequalities extend to
C$^*$-algebras.

\begin{theorem} If $A$ is a C$^*$-algebra, if $\tau$ is a faithful
tracial
state on $A$, and if $p$ and $q$ are positive real numbers for
which $p^{-1}+q^{-1}=1$,
then, for all $a,b\in A^+$,
\begin{enumerate}
\item $\tau(|ab|)\le p^{-1}\tau(a^p)\,+\,q^{-1}\tau(b^q)$, where equality holds
if and only if $b^q=a^p$, and
\item $\tau(|ab|^{\frac{1}{2}})\le\sqrt{\tau(a)\tau(b)}
\le\frac{1}{2}(\,\tau(a)+\tau(b)\,)$.
\end{enumerate}
\end{theorem}

\noindent{\em Proof.} Let $\tau$ be a faithful trace on $A$ and write
$\tau$ in its GNS representation $\tau(x)=\langle\pi_{\tau}(x)\xi_{\tau},
\xi_{\tau}\rangle$, for $x\in A$, where $\pi_{\tau}:A\rightarrow B(H_{\tau})$
is a $*$-representation of $A$ on a Hilbert space $H_{\tau}$ and $\xi_{\tau}
\in H_{\tau}$ is a unit cyclic vector for $\pi_{\tau}(A)$. If
$M=\pi_{\tau}(A)^{''}$
(the double commutant), then there is a faithful normal trace $\t(\cdot)$
on $M$ such that $\t(\,\pi_{\tau}(x)\,)=\tau(x)$, for all $x\in A$
(Proposition V.3.19 of \cite{takesaki}). Thus,
the results concerning tracial Young-type inequalities
in semifinite von Neumann algebras
apply to $A$ and $\tau$ through $\t(\cdot)$.
\qed\medskip

\section*{Acknowledgement}\label{S:acknowledgement}

We wish to thank Mart\'\i n Argerami and Juliana Erlijman
for their advice and
suggestions regarding the results herein.


\end{document}